\documentclass[lettersize,journal]{IEEEtran}
\usepackage{amsmath,amsfonts}
\usepackage{algorithm}
\usepackage{array}
\usepackage[caption=false,font=normalsize,labelfont=sf,textfont=sf]{subfig}
\usepackage{textcomp}
\usepackage[noend]{algpseudocode}
\usepackage{color}

\usepackage{amsmath}
\usepackage{stfloats}
\usepackage{url}
\usepackage{multirow}
\usepackage{verbatim}
\usepackage{graphicx}
\usepackage{cite}
\usepackage{booktabs}
\usepackage{mathtools}
\begin{document}
\title{High-Order Coupled Fully-Connected Tensor\\Network Decomposition for Hyperspectral\\ Image Super-Resolution}
\author{Diyi Jin, Jianjun  Liu, ~\IEEEmembership{Member~IEEE,} Jinlong Yang, Zebin Wu, \IEEEmembership{~Senior Member IEEE}
}
\maketitle
\begin{abstract}
Hyperspectral image super-resolution addresses the problem of fusing a low-resolution hyperspectral image (LR-HSI) and a high-resolution multispectral image (HR-MSI) to produce a high-resolution hyperspectral image (HR-HSI). Tensor analysis has been proven to be an efficient method for hyperspectral image processing. However, the existing tensor-based methods of hyperspectral image super-resolution like the tensor train and tensor ring decomposition only establish an operation between adjacent two factors and are highly sensitive to the permutation of tensor modes, leading to an inadequate and inflexible representation. In this paper, we propose a novel method for hyperspectral image super-resolution by utilizing the specific properties of high-order tensors in fully-connected tensor network decomposition. The proposed method first tensorizes the target HR-HSI into a high-order tensor that has multiscale spatial structures. Then, a coupled fully-connected tensor network decomposition model is proposed to fuse the corresponding high-order tensors of LR-HSI and HR-MSI. Moreover, a weighted-graph regularization is imposed on the spectral core tensors to preserve spectral information. In the proposed model, the superiorities of the fully-connected tensor network decomposition lie in the outstanding capability for characterizing adequately the intrinsic correlations between any two modes of tensors and the essential invariance for transposition. Experimental results on three data sets show the effectiveness of the proposed approach as compared to other hyperspectral image super-resolution methods.
\end{abstract}
\section{intorduction}
\IEEEPARstart{H}yperspectral images(HSIs) contain hundreds of spectral bands which range from visible to infrared wavelengths. With abundant spectral information, HSIs are used to many common applications, such as classification\cite{camps2013advances} and object detecion\cite{stein2002anomaly}. However, due to actual hardware constraints, it is hard to obtain HSI with high-spatial and high-spectral resolution simultaneously. A feasible solution is to fuse available low-resolution HSI (LR-HSI) and the panchromatic or multispectral images (HR-MSI) over the same scene, called HSI super-resolution (HSI-SR).

There are many types of researches in related fusion problems\cite{dian2021recent} \cite{yokoya2017hyperspectral} \cite{loncan2015hyperspectral}. Generally, they can be divided into four classes: component substitution\cite{tu2001new}, multi-resolution \cite{nencini2007remote}, model-based approaches\cite{dian2019nonlocal}, and learning-based approaches\cite{liu2022model}. Many approaches focus on the last two classes. For the model-based approaches, they address the HSI-SR problem by optimizing the model based on various priors and the degradation relationships consisting of matrix factorization methods \cite{liu2020truncated}\cite{xue2021spatial}, tensor-based methods, and so on. For the learning-based approaches, they are data-driven and can learn the image features from the training data.

Recently, tensor-based methods have achieved significant performance in many applications\cite{xu2020hyperspectral}\cite{xue2022laplacian}\cite{xue2021multilayer}. In the HSI-SR problem, they have advantages preserving the spatial-spectral structure of HR-HSI. The common tensor decomposition methods are canonical polyadic decomposition \cite{kanatsoulis2018hyperspectral}, Tucker decomposition \cite{li2018fusing}, the tensor train decomposition \cite{dian2019learning} and the tensor ring decomposition \cite{chen2021hyperspectral}. However, these methods still have some weaknesses. On the one hand, some deep structural information like multiscale correlation is not considered. On the other hand, these decompositions only establish a connection between adjacent two factors, and thus are highly sensitive to the permutation of tensor modes.

In this paper, to capture the underlying structure of the HR-HSI, we propose a novel model for HSI-SR by utilizing the specific properties of high-order tensor via fully-connected tensor network decomposition (FCTN)\cite{zheng2021fully}. Firstly, we unfold the LR-HSI and HR-MSI into high-order tensors which can describe the multiscale correlation of patches. Secondly, the FCTN model is proposed to fuse the LR-HSI and HR-MSI. Finally, we impose a weighted graph regularization (WGR) on the spectral mode to maintain the spectral structure of the target HR-HSI. Compared to other existing approaches, the method has innovative characteristics as follows:

1) We unfold the LR-HSI and HR-MSI into high-order tensors that can describe the intrinsic information between different scales in both spatial and spectral dimensions.

2) The FCTN can make better use of the correlations between different factors and reveal the dimensional correlations more completely.

3) To maintain the spectral structure, we propose a weighted-graph regularization to constrain the model.
\section{Proposed Approach}
\subsection{High-Order Tensorization}
In this article, the target HR-HSI is denoted by $\mathcal{X}$ $\in{\mathbb{R}}$ $^{M\times N\times S}$, where $S$ represents the number of bands, ${M}$ and ${N}$ represent the number of rows and columns, respectively. Correspondingly,  $\mathcal{Y}$ $\in\mathbb{R}^{\mathit{m}\times {n}\times{S}}$ denotes the LR-HSI, which is a spatially degraded image $\mathcal{X}$ and satisfies $m=M/p$, $n=N/p$ with $p$ being the downsamping ratio. $\mathcal{Z}$ $\in \mathbb{R}^{\mathit{M}\times {N}\times{s}}$ denotes the HR-MSI, which is spectrally downsampled concerning $\mathcal{X}$, satisfying $S>s$. The goal of fusion is to estimate $\mathcal{X}$ from the observations $\mathcal{Y}$ and $\mathcal{Z}$.

Both $\mathcal{Z}$ and $\mathcal{Y}$ are the downsampled versions of $\mathcal{X}$, so they can be represented as $\mathcal{Z}_{(3)}={\mathbf R}\mathcal{X}_{(3)}$ and $\mathcal{Y}_{(3)}=\mathcal{X}_{(3)}{\mathbf P}$. Here $\mathcal{X}_{(3)}\in\mathbb{R}^{S\times{MN}}$ represents the 3-mode of tensor $\mathcal{X}$ by unfolding the tensor into a matrix along the third mode. $\mathbf{R}\in\mathbb{R}^{s\times S}$ and $\mathbf{P}\in\mathbb{R}^{MN\times mn}$ represent the spectral response function (SRF) of sensor and spatial down-sampling operator, respectively.

The high-order tensorization of $\mathcal{X}$ consists of two steps. In the first step, spatial dimensions are divided into $M=M_1\times M_2\times\cdots\times M_d$ and $N=N_1\times N_2\times\cdots\times N_d$, and then make the size ${M\times N\times S}$ changed into ${M_1\times\cdot \cdot \cdot M_d\times N_1\times\cdots{N_d}\times S}$. In the second step, we use permution operation to obtain the high-order tensor of size ${{M_1N_1}\times{M_2N_2}\times\cdot\cdot\cdot\times{M_dN_d} \times S}$ following the column-first rules. It is shown that each mode of the high-order tensor represents the patches at a different scale. For the convenience of presentation, we use $\mathcal{T}\left\{\right\}$ to represent the tensorization steps, so the high-order tensor is $\mathcal{T}\left\{\mathcal{X}\right \} \in \mathbb{R}^{M_1N_1\times M_2N_2\times\cdot\cdot\cdot\times M_dN_d\times S}$.

Analogously, $\mathcal{Y}$ $\in\mathbb{R} ^{\mathit{m}\times {n}\times{S}}$ can be represented as $\mathcal{T}\left\{\mathcal{Y}\right\}\in\mathbb{R}^{m_1n_1\times m_2n_2\times\cdot\cdot\cdot\times m_dn_d \times S}$, which is spatially downsampled $\mathcal{T}\left\{\mathcal{X}\right\}$. $\mathcal{Z}$ $\in \mathbb{R}^{\mathit{M}\times {N}\times{s}}$ can be represented as $\mathcal{T}\left\{\mathcal{Z}\right \} \in \mathbb{R}^{M_1N_1\times M_2N_2\times\cdot\cdot\cdot\times M_dN_d\times s}$, which is spectrally downsampled $\mathcal{T}\left\{\mathcal{X}\right\}$.
\subsection{FCTN Representation}
The FCTN representation aims to decompose $\mathcal{T}\left\{\mathcal{X}\right\}\in\mathbb{R}^{M_1N_1\times{M_2N_2}\times\cdots\times{S}}$ into a set of $\left(d+1\right)$-mode factor tensors $\mathcal U_t\in\mathbb{R}^{\prod_{i=1}^{t-1}r_{t,i}\times{M_tN_t}\times\prod_{j=t+1}^{d+1}r_{t,j}}$, $t\in\left\{{1,2,\cdots,d+1}\right\}$, where $r_{t,t+1}$ represents the rank between $\mathcal U_t$ {and} $\mathcal U_{t+1}$. Let $\mathcal{T}\left\{\mathcal{X}\right\}{(i_1,i_2,\cdots,i_{d+1})}$ be the $(i_1,i_2,\cdots,i_{d+1})$th element of $\mathcal{T}\left\{\mathcal{X}\right\}$. The FCTN representation\cite{zheng2021fully} to represent $\mathcal{T}\left\{\mathcal{X}\right\}$ is defined as follows:
\begin{equation}
\begin{array}{l}
\mathcal{T}\left\{\mathcal{X}\right\}{(i_1,i_2,\cdots,i_{d+1})}=\\
\sum_{j_{1,2}=1}^{r_{1,2}}\sum_{j_{1,3}=1}^{r_{1,3}}\cdots\sum_{j_{1, d+1}=1}^{r_{1,d+1}}\sum_{j_{2,3}=1}^{r_{2,3}}\cdots\sum_{j_{d,d+1}=1}^{r_{d, d+1}}\\\mathcal{U}_{1}\left(i_1, j_{1,2},j_{1,3},\cdots, j_{1, d+1}\right)\\\mathcal{U}_{2}\left(j_{1,2},i_2,j_{2,3}, \cdots,j_{2, d+1}\right)\cdots\\\mathcal{U}_{d+1}\left(j_{1,d+1},j_{2, d+1},\cdots,j_{d, d+1},i_{d+1}\right)
\label{rff}
\end{array}
\end{equation}

Conveniently, denoting FCTN representation by $\mathcal{F}()$, we obtain $\mathcal{T}\left\{\mathcal{X}\right \}=\mathcal{F}(\mathcal U_1,\mathcal U_2,\cdots,\mathcal U_{d+1})+{\varepsilon}$ where $\varepsilon$ is the error. Similarly, $\mathcal{T}\left\{\mathcal{Z}\right\}$ and $\mathcal{T}\left\{\mathcal{Y}\right\}$ can be represented as $\mathcal{T}\left\{\mathcal{Z}\right\}=\mathcal{F}(\mathcal U_1,\mathcal U_2,\cdots,\mathcal U_{d+1}\times_{d+1}\mathbf{R})+{\varepsilon}$ and $\mathcal{T}\left\{\mathcal{Y}\right\}=\mathcal{F}(\mathcal Q,\mathcal U_2,\cdots,\mathcal U_{d+1})+{\Large \varepsilon}$ where $\mathcal Q$ represents the downsampling of $\mathcal U_1$ and $\times_{d+1}$ denotes the multiplication between tensor and matrix in the ${(d+1)}$th-mode of tensor. And we can formulate the HSI-SR problem as follows:
\begin{equation}
\begin{aligned}
\min_{\mathcal{Q},\mathcal U_{t,t=1...d+1}}\frac{1}{2}||\mathcal{T}\left \{\mathcal{Y}\right \}-\mathcal{F}(\mathcal{Q},\mathcal U_2,\cdots,\mathcal U_{d+1})||_F^2\\+\frac{\lambda}{2}||\mathcal{T}\left\{\mathcal{Z}\right \}-\mathcal{F}(\mathcal U_1,\mathcal U_2,\cdots,\mathcal U_{d+1}\times_{d+1}\mathbf{R})||_F^2
\label{formula_8}
\end{aligned}
\end{equation}
where the $\mathit{\lambda}$ represents the parameter to balance of two terms and $\mathcal{||\mathcal{X}|}|_F=\left(\sum{_{i_1,i_2,i_3\cdots{i_{d+1}}}}\mathcal{X}\left(i_1,i_2,i_3\cdots{i_{d+1}}\right)^2\right)^{1/2}$ is Frobenius norm.
\subsection{Weighed-Graph Regularization}
The core tensor $\mathcal{U}_{d+1}$ is related to the spectral structure of $\mathcal{X}$. To keep the spectral structure of $\mathcal{X}$, the WGR is imposed on the $\mathcal{U}_{d+1}$. The graph is built as $\mathcal{G}$=$(\mathbf{V},\mathbf{E},\mathbf{W})$, where $\mathbf{V}$ is the set of vertices and $\mathbf{E}$ is the set of edges standing for the bands of $\mathcal{Y}$. $\mathbf{W}$ is the set of weights in measuring the closeness of two bands and defined as follows:
\begin{equation}
\mathbf{W}{\left(i,j\right)}=\left\{\begin{array}{ll}
\exp \left(-\frac{\|\mathcal{Y}(:,:, i)-\mathcal{Y}(:,:, j)\|_{F}^{2}}{\sigma^{2}}\right), & j \in\Omega\\
0, &\text{otherwise}
\end{array}\right.
\end{equation}
where $\sigma$ controls the degree of smoothness and $\mathit{j}\in\Omega$ represents the adjacent of band $\mathit{i}$. The WGR can be written as follows:
\begin{equation}
\begin{aligned}
&\frac{1}{2}\sum_{ij}\|\mathcal U_{d+1(d+1)}(i,:)-\mathcal U_{d+1(d+1)}(j,:) \|_{F}^{2}\mathbf{W}{\left(i,j\right)}\\
&={\mathrm{Tr}}(\mathcal U_{{d+1}{(d+1)}}^T\mathbf{L_S}\mathcal U_{{d+1}{(d+1)}})
\end{aligned}
\end{equation}
where $\mathbf{L_S}$ is Laplacian matrix computed by $\mathbf{D-W}$, $\mathbf{D}=diag(\mathbf{I_1,I_2,...,I_S})$ and $\mathbf{I_i}={\sum_j\mathbf{W}}{\left(i,j\right)}$. By integrating the WGR into (\ref{formula_8}) and it can be written as follows:
\begin{equation}
\begin{aligned}
&\min_{\mathcal Q,\mathcal U_{t,t=1...d+1}}\frac{1}{2} || \mathcal{T}\left \{ \mathcal{Y}  \right \}-\mathcal{F}(\mathcal Q,\mathcal U_2,...,\mathcal U_{d+1})||_F^2 \\
&+\frac{\lambda}{2}||\mathcal{T}\left \{ \mathcal{Z}\right\}-\mathcal{F}(\mathcal U_1,\mathcal U_2,...,\mathcal U_{d+1}\times_{d+1}\mathbf{R})||_F^2 \\
&+\frac{\beta }{2}\rm Tr(\mathcal U_{{d+1}{(d+1)}}^T\mathbf{L_S}\mathcal U_{{d+1}{(d+1)}})+\frac{\mu}{2}\left(\sum_{t=1}^{d}||\mathcal U_t||_F^2+||\mathcal Q||_F^2\right)
\label{formula_11}
\end{aligned}
\end{equation}
where ${\beta}$ and ${\mu}$ are the parameters that control the importance of the WGR and tensor decay regularization, respectively.
\subsection{Optimization Algorithm}
To solve (\ref{formula_11}), the alternating optimization framework can be adopted. According to the FCTN representation, if one of $\mathcal U_t, {t\in\left\{{1,2,\cdots,{d+1}}\right\}}$ does not participate in the composition of $\mathcal{F}\left\{\mathcal{Z}\right\}$, we denote it by $\mathcal O^{\neq{t}}=\mathcal{F}(\mathcal U_1,\mathcal U_2,\cdots,\mathcal U_{t-1},\mathcal U_{t+1},\cdots,\mathcal U_{d+1}\times_{d+1}\mathbf{R})$, and if $\mathcal U_t$ does not participate in the composition of $\mathcal{F}\left\{\mathcal{Y}\right\}$, we denote it by $\mathcal H^{\neq t}=\mathcal{F}(\mathcal Q,\cdots,\mathcal U_{t-1},\mathcal U_{t+1}, \cdots,\mathcal U_{d+1})$. Then, we can gain the relation as follows:
\begin{equation*}
\mathcal{T}\left\{\mathcal{Z}\right\}_{(t)}=\mathcal{U}_{t(t)}\left(\mathcal{O}^{\neq t}_{\left(b_{1:d};c_{1:d}\right)}\right)^T
\end{equation*}
\begin{equation*}
\mathcal{T}\left\{\mathcal{Y}\right\}_{(t)}=\mathcal {U}_{t(t)}\left(\mathcal {H}^{\neq t}_{\left(b_{1:d};c_{1:d}\right)}\right)^T
\end{equation*}
\begin{equation*}
b_{i}=\{\begin{array}{ll}2 i, & \text { if } i<t, \\
2 i-1, & \text { if } i \geq t,
\end{array} \text { and } c_{i}=\{\begin{array}{ll}2 i-1, & \text { if } i<t \\
2 i, & \text { if }i\geq t
\end{array}
\end{equation*}
Here $\mathcal{O}^{\neq{t}}_{\left(b_{1:d};c_{1:d}\right)}$ is gained by reshaping the size of $\mathcal O^{\neq t}$ into multiplication of $\prod_{i=1,\neq{t}}^{d}M_{b_i}N_{b_i}$ and $\prod_{j=1,\neq{t}}^{d}r_{t,c_j}$. Analogously, $\mathcal{H}^{\neq t}_{\left(b_{1:d};c_{1:d}\right)}$ have the similar definition. Then, we can optimize the proposed model by solving the following subproblems.\\
1) Solving $\mathcal{Q}$ subproblem: Optimizing with respect to $\mathcal{Q}$ can be written as
\begin{equation}
\min_{\mathcal{Q}_{(1)}}\frac{\lambda}{2}\left\|\mathcal{T}\{\mathcal{Y}\}_{(1)}-\mathcal{Q}_{(1)}(\mathcal{H}_{{\left(b_{1:d};c_{1:d}\right)}}^{\neq 1})^T\right\|_{F}^{2}+\frac{\mu}{2}\left\|\mathcal{Q}_{(1)}\right\|_{F}^{2}
\end{equation}
the solution of above is
\begin{equation}
\begin{array}{l}
\mathcal{Q}_{(1)}=\left(\lambda \mathcal{T}\{\mathcal{Y}\}_{(1)} \mathcal{H}_{\left(b_{1: d} ; c_{1: d}\right)}^{\neq 1}\right)\left(\lambda\left(\mathcal{H}_{\left(b_{1: d} ; c_{1: d}\right)}^{\neq 1}\right)^{T}\right. \\
\left.\mathcal{H}_{\left(b_{1: d} ; c_{1: d}\right)}^{\neq 1}+\mu I\right)^{-1}
\end{array}
\end{equation}
2) Solving $\mathcal U_1$ subproblem: Optimizing with respect to $\mathcal U_1$ is written as
\begin{equation}
\min_{\mathcal{U}_{1(1)}}\frac{1}{2}\|\mathcal{T}\{\mathcal{Z}\}_{(1)}-\mathcal{U}_{1(1)}\left(\mathcal{O}_{\left(b_{1:d};c_{1:d}\right)}^{\neq 1}\right)^{T}\|_{F}^{2}+\frac{\mu}{2}\|\mathcal{U}_{1(1)}\|_{F}^{2}
\label{lll}
\end{equation}
the solution of (\ref{lll}) is calculated by
\begin{equation}
\begin{array}{l}
\mathcal{U}_{1(1)}=\left(\lambda \mathcal{T}\{\mathcal{Z}\}_{(1)} \mathcal{O}_{\left(b_{1: d} ; c_{1: d}\right)}^{\neq 1}\right)\left(\lambda\left(\mathcal{O}_{\left(b_{1: d} ; c_{1: d}\right)}^{\neq 1}\right)^{T}\right. \\
\left.\mathcal{O}_{\left(b_{1: d} ; c_{1: d}\right)}^{\neq 1}+\mu I\right)^{-1}
\end{array}
\end{equation}
3) Solving $\mathcal U_t$, $\mathit{t}$=2...$\mathit{d}$, subproblem: Optimizing with respect to $\mathcal U_t$, $\mathit{t}$=2...$\mathit{d}$ is written as
\begin{equation}
\begin{aligned}
&\min_{\mathcal{U}_{t}}\frac{1}{2} \left\|\mathcal{T}\{\mathcal{Z}\}_{(t)}-\mathcal{U}_{t(t)}(\mathcal{O}_{{\left(b_{1:d}; c_{1:d}\right)}}^{\neq t})^{T}\right\|_{F}^{2}\\
&+\frac{\lambda}{2}\left\|\mathcal{T}\{\mathcal{Y}\}_{(t)}-\mathcal{U}_{t(t)}(\mathcal{H}_{{\left(b_{1:d};c_{1:d}\right)}}^{\neq t})^{T}\right\|_{F}^{2}+\frac{\mu}{2}\left\|\mathcal{U}_{t(t)}\right\|_{F}^{2}
\label{ttt}
\end{aligned}
\end{equation}
the solution of (\ref{ttt}) is calculated by:
\begin{equation}
\begin{aligned}
\mathcal{U}_{t(t)}=&\left(\mathcal{T}\{\mathcal{Z}\}_{(t)} \mathcal{O}_{\left(b_{1: d} ; c_{1: d}\right)}^{\neq t}+\lambda \mathcal{T}\{\mathcal{Y}\}_{(t)} \mathcal{H}_{\left(b_{1: d} ; c_{1: d}\right)}^{\neq t}\right) \\
&\left(\left(\mathcal{O}_{\left(b_{1: d} ; c_{1: d}\right)}^{\neq t}\right)^{T} \mathcal{O}_{\left(b_{1: d} ; c_{1: d}\right)}^{\neq t}+\lambda\left(\mathcal{H}_{\left(b_{1: d} ; c_{1: d}\right)}^{\neq t}\right)^{T}\right.\\
&\left.\mathcal{H}_{\left(b_{1: d} ; c_{1: d}\right)}^{\neq t}+\mu I\right)^{-1}
\end{aligned}
\end{equation}
4) Solving $\mathcal U_{(d+1)}$ subproblem: Optimizing with respect to $\mathcal U_{(d+1)}$ is written as
\begin{equation}
\begin{aligned}
\min_{\mathcal U_{d+1}}\frac{1}{2} &\left\|\mathcal{T}\{\mathcal{Z}\}_{(d+1)}-\mathbf{R} \mathcal{U}_{d+1(d+1)}(\mathcal{O}_{\left(b_{1:d};c_{1:d}\right)}^{\neq d+1})^{T}\right\|_{F}^{2} \\
&+\frac{\lambda}{2}\left\|\mathcal{T}\{\mathcal{Y}\}_{(d+1)}-\mathcal{U}_{d+1(d+1)}(\mathcal{H}_{\left(b_{1:d};c_{1:d}\right)}^{\neq d+1})^{T}\right\|_{F}^{2} \\
&+\frac{\beta}{2}\operatorname{Tr}\left(\mathcal{U}_{d+1(d+1)}^{T}\mathbf{L_{S}}\mathcal{U}_{d+1(d+1)}\right)
\label{yyy}
\end{aligned}
\end{equation}
By making the gradient zero, we have the following equation:
\begin{equation}
\begin{aligned}
&\mathbf{R}^{T}\mathbf{R}\mathcal{U}_{{d+1}_{(d+1)}} \left(\mathcal{O}_{{\left(b_{1:d};c_{1:d}\right)}}^{\neq d+1}\right)^T\mathcal{O}_{{\left(b_{1:d};c_{1:d}\right)}}^{\neq d+1}+\\&\lambda\mathcal{U}_{d+1(d+1)}\left(\mathcal{H}_{{\left(b_{1:d};c_{1:d}\right)}}^{\neq{d+1}}\right)^T\mathcal{H}_{{\left(b_{1:d};c_{1:d}\right)}}^{\neq{d+1}}+\beta\mathbf{L_S}\mathcal{U}_{d+1(d+1)}=\\&\mathbf{R}^{T}\mathcal{T}\{\mathcal{Z}\}_{(d+1)}\mathcal{O}_{{\left(b_{1:d};c_{1:d}\right)}}^{\neq{d+1}}+\lambda\mathcal{T}\{\mathcal{Y}\}_{(d+1)}\mathcal{H}_{{\left(b_{1:d};c_{1:d}\right)}}^{\neq d+1}
\end{aligned}
\end{equation}
The equation can be solved by using conjugate gradient method.\\
\begin{table}
\begin{tabular}{ll}
\hline
 {\bf{Algorithm 1}} FCTN for HSI-SR\\
\hline
1:{\bf{Input:}} HR-MSI $\mathcal{Z}$, LR-HSI $\mathcal{Y}$,
parameters $\mathit{\mu}$, $\mathit{\lambda}$, $\mathit{\beta}$, $\sigma$, $M_t, N_t$, $r_{t,t}$ \\
2:{\bf{Initialization:}}\\
$\mathcal U_{t}=rand(r_{t,1},r_{t,2}...r_{t,t-1},M_tN_t,r_{t,t+1},..r_{t,{d+1}})$.\\
$Q=rand({{m_1n_1}},r_{1,2},r_{1,3}.. r_{1,d},r_{1,{d+1}})$.\\
$\mathcal{T}\{{Z}\}$=$\mathcal{F}$(${U}_{1},{U}_{2},\dots,{U}_{d+1}\times_{d+1}\bf{R}$).\\
$\mathcal{T}\{{Y}\}$=$\mathcal{F}$($Q,{U}_{2},\ldots,{ U_{d+1}}$).\\
3:Computing the weight matrix ${\bf{W}}$ by (3).\\
4:{\bf{while}} $Iteration$ $\le maxiter$ \\
5:\quad{Update} $\mathcal{Q}$ by (7).\\
6:\quad{Update} $\mathcal U_1$ by (9).\\
7:\quad\bf{for} $\mathit{t}=2:{d}$\\
8:\quad\quad{ Update} $\mathcal U_t$ by (11).\\
9:\quad{\bf{end for}} \\
10:\quad{Update} ${\mathcal U_{d+1}}$ by (13)\\
11:\quad$Iteration=Iteration$+1.\\
12:{ \bf{end while}}\\
13:{\bf{Output}}:$\mathcal{X}$ = $\mathcal{T}^{-1}$ $\left \{\mathcal{F}( U_{1}, U_{2},..., U_{d+1}) \right\}$.\\
\hline
\end{tabular}
\end{table}
After getting the latent core tensors, we apply the $\mathcal{T}^{-1}\left\{\right\}$ operation to rebound the target $\mathcal{X}$, where $\mathcal{T}^{-1}\left\{\right\}$ represents the reverse reconstruction of $\mathcal{T}\left\{\right\}$. Details of the algorithm are given in Algorithm 1.
\section{Experiments}
\subsection{Data Sets}
The first data set is the Chikusei data set captured by HeadWall’s Hyperspec Visible and Near-Infrared, series C imaging sensor over Chikusei in Japan. The size of image is 2517$\times$2355$\times$128 covering the spectral range from 0.363 to 1.018 $\mu$m. The sub-image with a size of 240$\times$240$\times$128 is cropped as the reference image.

The second data set is the SanDiego data set which generated by the Airborne Visible/Infrared Imaging Spectrometer (AVIRIS) sensor in San Diego, CA, USA. After removing low SNR bands and water absorption bands, the number of the spectral bands is 186. In the experiments, we select the top-left part of size 200$\times$200$\times$186 as the reference image.

The last data set is the University of Pavia which was gained by the Reflective Optics System Imaging Spectrometer (ROSIS) in Pavia, Italy. The size of the image is 610$\times$340$\times$115 covering the spectrums from 0.43 to 0.86 $\mu$m. It contains 103 bands after removing the water vapor absorption bands. We select a sub-image with a size of 256$\times$256$\times$103 for the experiments.

To generate the LR-HSI, we spatially blur the reference image by averaging 8$\times$8 disjoint spatial blocks and then downsample the blurred image by a factor of 8 in two spatial directions. Considering the different wavelength range and number of bands, we use different SRFs in different data sets. In the Chikusei data set, we generate the SRF by referring to Hysure\cite{simoes2014convex} method and obtain an eight-band HR-MSI. In the SanDiego data set, we generate the HR-MSI by averaging the bands of the reference image according to the SRF of IKONOS. In the University of Pavia data set, we use an IKONOS-like SRF to generate a four-band HR-MSI. In order to simulate the real fusion as much as possible, Gaussian noise is added to all the generated LR-HSI (SNR=25dB) and HR-MSI (SNR=25dB).
\begin{figure}[htpb]
\small
\centering
\subfloat[]{\includegraphics[width=1.6in]{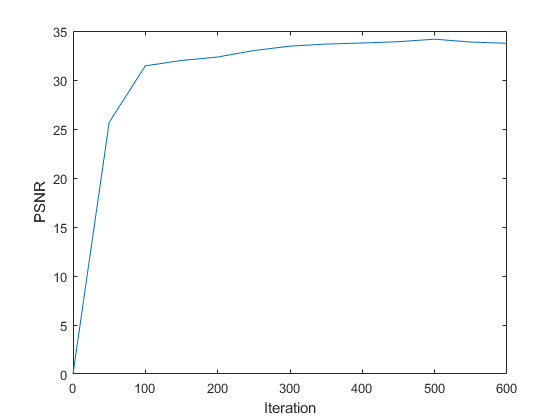}}
\hfil
\subfloat[]{\includegraphics[width=1.6in]{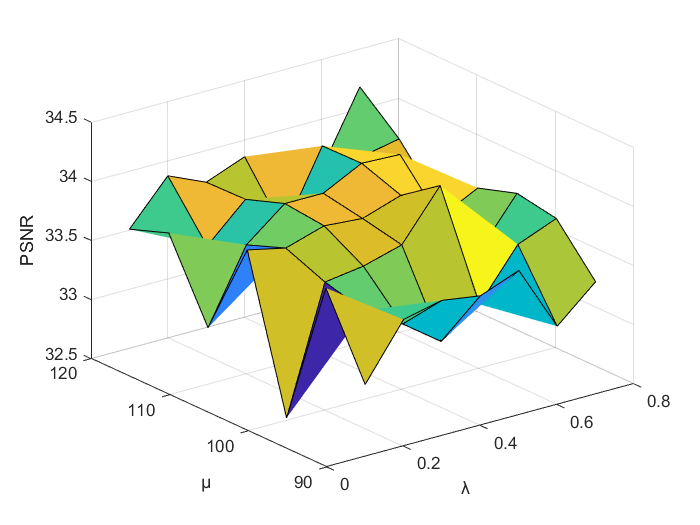}}
\\
\caption{Results for SanDiego data set. (a) The PSNR curve as a function of ${\mathit{Iteration}}$. (b) The PSNR curve as a function of $\mathit{\mu}$ and $\mathit{\lambda}$.
}
\label{fig_3}
\end{figure}
\subsection{Setting of Parameters and Quantitative Metrics}
In this section, we investigate effect of the parameters in the proposed method on SanDiego data set. Fig.1 illustrates the PSNR results of proposed method as a function of $Iteration$ and a function of $\mathit{\lambda}$, $\mathit{\mu}$. For the proposed FCTN, we set $\lambda$=0.1, $\mu$=120, $Iteration$=480, $\sigma=10$, $\beta=0.1$. In our experiment, there is a widely optimal range for the choice of orders, and ranks. For the University of Pavia data set, the dimension of each order is $M_1=N_1=8,M_2=N_2=8,M_3=N_3=2,M_4=N_4=2$ to [64$\times$64$\times$4$\times$4$\times$103]. For the rank, we set $r_{1,2}$ from 36 to 40, $r_{1,3}=r_{1,4}$=3, $r_{i,5}$=2 with $i$ from 1 to 5, $r_{2,3}$=9 and set $r_{2,4}=r_{3,4}$=4. For the Chikusei data set, the dimension of each order is $M_1=N_1=8,M_2=N_2=5,M_3=N_3=2, M_4=N_4=3$ to [64$\times$25$\times$4$\times$9$\times$128]. For the rank, we set $r_{1,2}$ from 36 to 40, $r_{1,3}=r_{1,4}$=3, $r_{i,5}$=2 where $i$ from 1 to 5, $r_{2,3}$ is set to 9 and $r_{2,4}$=$r_{3,4}$=4. For the Sandiego data set, the dimension of each order is $M_1=N_1=8,M_2=N_2=5,M_3=N_3=5$ to [64$\times$25$\times$25$\times$186]. For the rank, we set $r_{1,2}$ from 42 to 46, $r_{1,3}=r_{1,4}$=2, $r_{2,3}$=8, $r_{2,4}=r_{3,4}$=5.

We use four common quantitative metrics to evaluate the quality of the reconstructed image. We bring into the peak-signal-to-noise (PSNR) to evaluate the spatial recovery quality in each band. The spectral angle mapper (SAM) is used to quantify spectral information preservation. The relatively dimensionless global error in synthesis (ERGAS) shows a global quality of the recovery image. The $Q2^{N}$ can jointly quantify spectral and spatial distortions.
\begin{table*}
\caption{QUALITY MEARSURE OF THREE DATASETS}
\centering
\resizebox{1.9\columnwidth}{!}{
\begin{tabular}{lcccccccccccc}
\hline
\multirow{2}{*}{\text{Method}}&
\multicolumn{4}{c}{\text{SANDIEGO}} &\multicolumn{4}{c}{UNIVERSITY OF PAVIA}&\multicolumn{4}{c}{CHIKUSEI}\\
\cline{2-13}&\text{PSNR}&\text{SAM} & \text{ERGAS} & \text{${Q{2^N}}$}&\text{PSNR}&\text{SAM} & \text{ERGAS} & \text{${Q{2^N}}$} &\text{PSNR}&\text{SAM} & \text{ERGAS} & \text{${Q{2^N}}$}\\
\hline
FCTN& \bf34.12& \bf2.91& \bf0.85 & \bf0.95 & 36.64& \bf4.26 & \bf1.24  & \bf0.80& \bf{36.38}& \bf{2.37}& \bf{1.70}&\bf{0.83}\\
CNMF& 31.97& 4.72& 1.08 & 0.92 & 34.41 & 5.42 & 1.52 & 0.77& 31.58&4.66&2.47& 0.68\\
LTMR& 30.77 & 4.29& 1.20 & 0.92 & 31.49 & 8.6  & 2.56  & 0.69 & 35.01& 2.9& 1.88&0.79\\
LTTR& 28.63& 4.52& 2.22  & 0.79& 32.15 & 4.13  & 2.56  & 0.67& 30.15& 6.3& 3.47& 0.75\\
HCTR& 30.52&3.18& 1.17& 0.92 &\bf36.70 & 4.45 & 1.34& 0.74 & 34.94& 3.25 & 2.48 & 0.79\\
Hysure& 28.99& 7.29& 2.13& 0.92 & 32.99 & 7.29 & 2.13  & 0.76& 28.99& 7.29& 2.13& 0.92\\
NPTSR& 32.50& 4.48&1.44& 0.92& 32.98 & 5.99 & 1.88  & 0.79 & 33.61& 3.01& 2.02& 0.88\\
\hline
\end{tabular}
}
\end{table*}
\begin{table}
\scriptsize
\centering
\caption{ABLATION STUDY ON USING WGR IN SANDIEGO DATASET}
\resizebox{.8\columnwidth}{!}{
\begin{tabular}{c|l|l|l|l}
\hline
WGR & \multicolumn{1}{c|}{PSNR} & \multicolumn{1}{c|}{SAM} & \multicolumn{1}{c|}{ERGAS} & \multicolumn{1}{c}{${Q{2^N}}$}\\
\hline
 w/ & \bf34.12& \bf2.91& \bf0.85 & \bf0.95\\
 w/o & 33.10& 3.20&  0.95 & 0.90\\
\hline
\end{tabular}}
\end{table}
\subsection{Results for Visual Images and Ablation Study}
\begin{figure*}[htpb]
\centering
\subfloat[]{\includegraphics[width=2in]{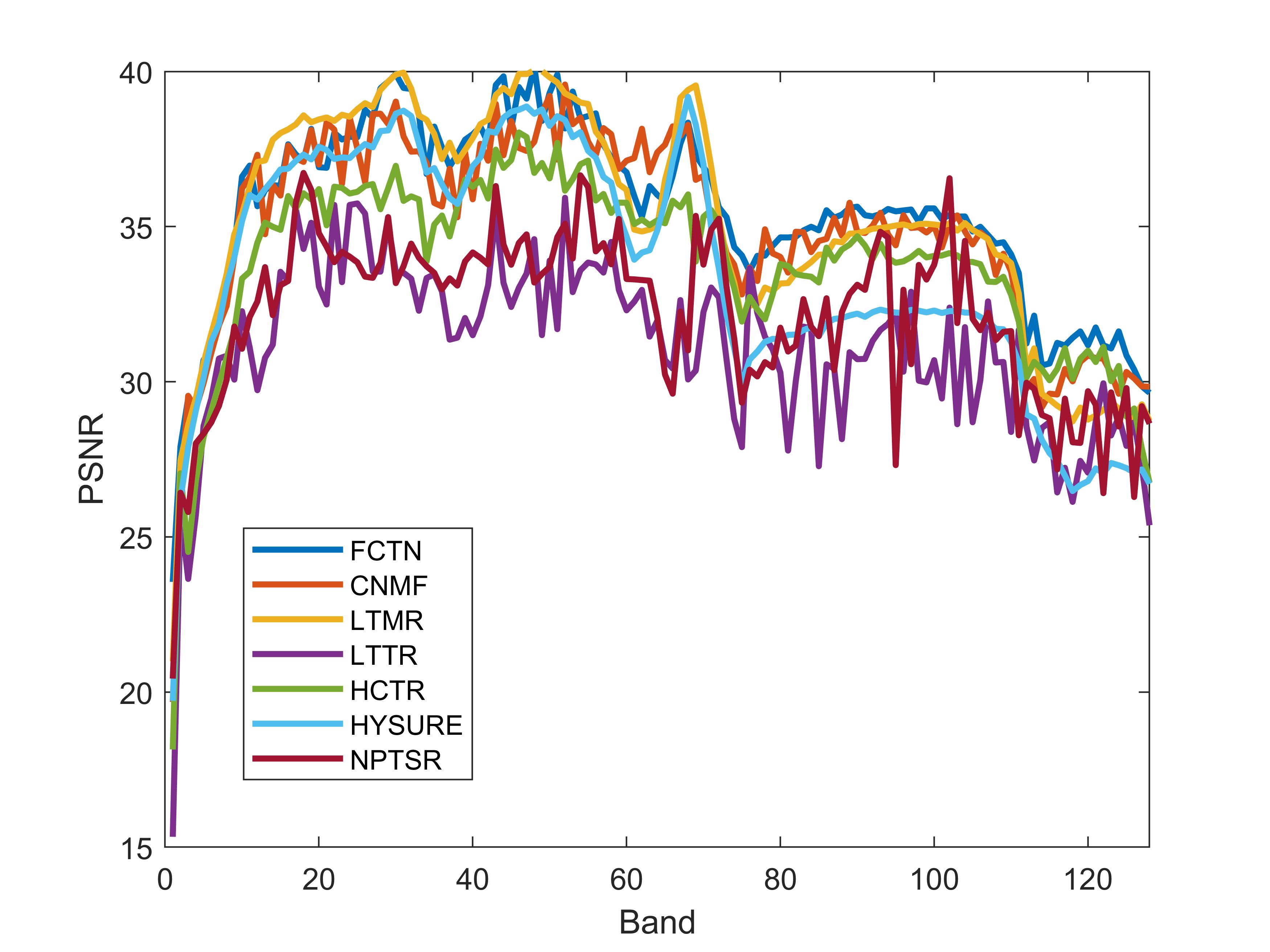}}
\hfil
\subfloat[]{\includegraphics[width=2in]{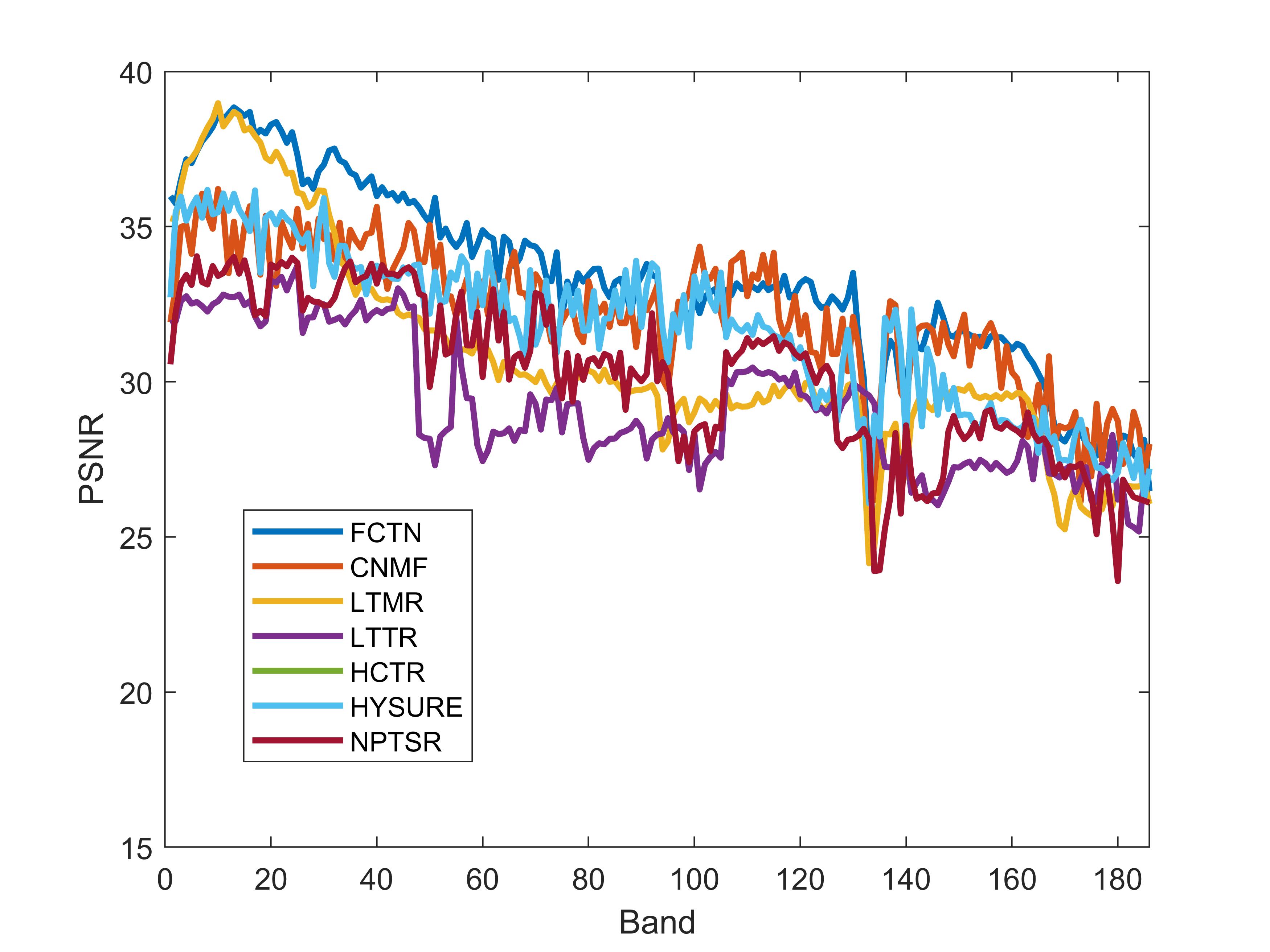}}
\hfil
\subfloat[]{\includegraphics[width=2in]{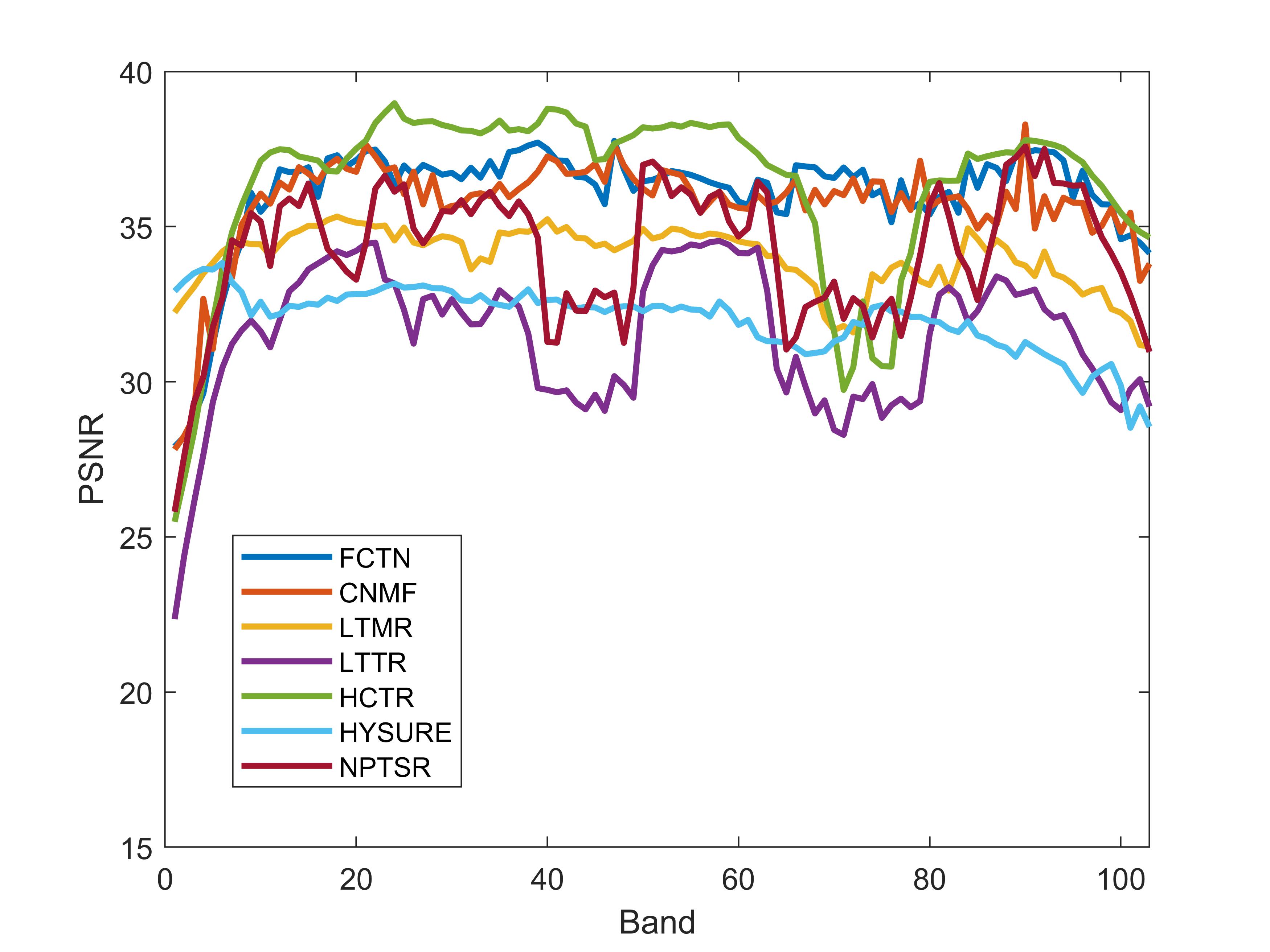}}
\caption{PSNR as a function of spectral band. (a) Chikusei data set. (b) SanDiego data set. (c) University of Pavia data set.}
\label{fig_1}
\end{figure*}

In this section, we evaluate the performance of proposed approach by comparing with six methods. The six methods are CNMF\cite{yokoya2011coupled}, LTMR \cite{dian2019hyperspectral}, LTTR\cite{dian2019learning}, HCTR \cite{xu2020hyperspectral}, Hysure\cite{simoes2014convex}, NPTSR\cite{xu2019nonlocal}. The experiment results are shown in Table I. Table shows the effects of all methods in four quantitative metrics and we marked best value as bold. Fig.\ref{fig_1} illustrates the fusion effect of the competing methods in different spectral bands. In the Chikusei, and SanDiego data sets, our proposed method achieved the best results. In the University of Pavia data set, our proposed method obtains the best results in the SAM, ERGAS, and $Q{2^N}$ quantitative metrics. The HCTR method is a relatively competitive method and gets the best effect in PSNR quantitative metric. Fig.\ref{fig_3} shows the reconstructed images and the residual images. Observing the residual image, we can find that the error distribution is relatively smooth, proving that our method has a good performance in restoraton.

The WGR is introduced to preserve the spectral structure. An ablation experiment is proposed to demonstrate the effectiveness of it. Its contribution  was assessed on the quality measures by removing it from (5). As shown in Table II, with  WGR performs better than without WGR, indicating that WGR contributes positively to the final result.
\section{Conclution}
In this paper, we propose the FCTN framework for HSI-SR. Compared to other methods, it is completely a novel method to deal with the high-order tensors. The FCTN framework showed its outstanding capability to adequately characterize the correlations between any two modes of tensors and was proved to be essentially transpositional invariable. Specifically, the FCTN framework consists of three steps. The first step is to transform the original data into high-order tensors. In the second step, we introduce the FCTN framework to solve the fusion problem of HR-MSI and LR-HSI. In the high-order tensor, some latent core tensors of HR-MSI and LR-HSI are shared. Finally, we use an alternating algorithm on FCTN framework. The results demonstrate the effectiveness of our proposed method.
\begin{figure*}[!t]
\centering
{\includegraphics[width=0.9in]{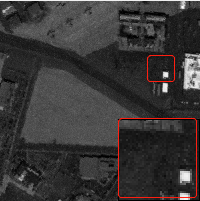}}
\hfil
{\includegraphics[width=0.9in]{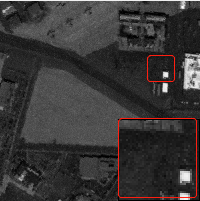}}
\hfil
{\includegraphics[width=0.9in]{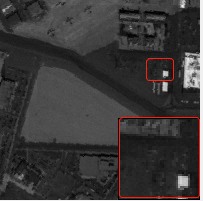}}
\hfil
{\includegraphics[width=0.9in]{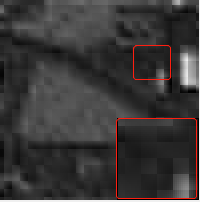}}
\hfil
{\includegraphics[width=0.9in]{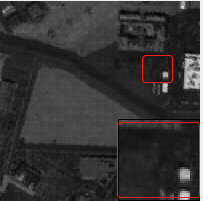}}
\hfil
{\includegraphics[width=0.9in]{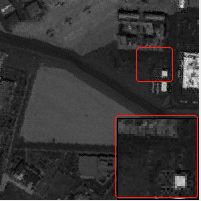}}
\hfil
{\includegraphics[width=0.9in]{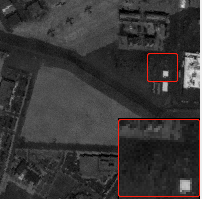}}
\\
\subfloat[]{\includegraphics[width=0.9in]{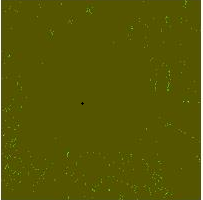}}
\hfil
\subfloat[]{\includegraphics[width=0.9in]{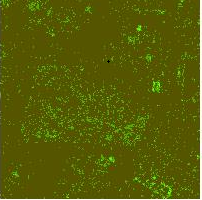}}
\hfil
\subfloat[]{\includegraphics[width=0.9in]{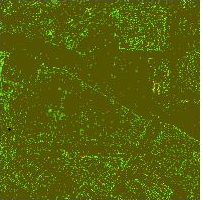}}
\hfil
\subfloat[]{\includegraphics[width=0.9in]{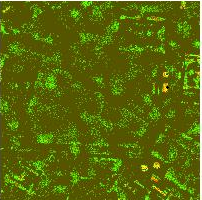}}
\hfil
\subfloat[]{\includegraphics[width=0.9in]{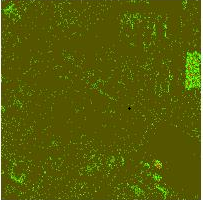}}
\hfil
\subfloat[]{\includegraphics[width=0.9in]{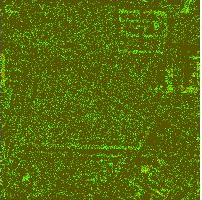}}
\hfil
\subfloat[]{\includegraphics[width=0.9in]{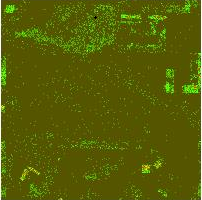}}
\\
\caption{Visual quality comparison for reconstructed images of the Sandiego data set. First row: reconstruct images at band 40th. Second row: residual images between the reference and reconstruct images. (a) FCTN. (b) CNMF. (c) LTMR. (d) LTTR. (e) HCTR. (f) Hysure. (g) NPTSR.}
\label{fig_3}
\end{figure*}
\bibliographystyle{IEEEtran}
\bibliography{my}
\end{document}